\documentclass[12pt,reqno]{amsart}
\usepackage{amsmath}
\usepackage{amsfonts}
\usepackage{amssymb}
\usepackage{caption}
\usepackage{latexsym}
\usepackage{color}
\usepackage{hyperref}
\usepackage{url}
\usepackage{enumerate}
\usepackage{booktabs}
\usepackage{longtable}
\usepackage[figuresleft]{rotating}

\usepackage{multirow}
\usepackage{lscape}
\usepackage{tabularx}

\renewcommand{\proof}{\noindent{\it Proof.\ \ }}
\renewcommand{\qed}{\ifmmode\square\else\nolinebreak\hfill
$\Box$\fi\par\vskip12pt}

\renewcommand\a{\alpha}

\newcommand\Ga{\mathrm{\Gamma}}   
\newcommand\Sig{{\it \Sigma}}

       \newcommand\D{\mathrm{D}}
    \newcommand\G{\mathrm{G}}   
\newcommand\K{\mathsf{K}}  \newcommand\M{\mathrm{M}}

\newcommand\ZZ{\mathbb{Z}}

          \newcommand\Aut{\mathrm{Aut}}
   \newcommand\Cay{\mathrm{Cay}}      \newcommand\Cos{\mathrm{Cos}}

\newcommand\soc{\mathrm{soc}}    \newcommand\Sy{\mathrm{S}}         
          
\newcommand\diam{\mathrm{diam}}            \newcommand\la{\langle}
\newcommand\ra{\rangle}
\newcommand\g{\mathbf{g}}

\newcommand\GL{\mathrm{GL}}                    
                    
\newcommand\PGL{\mathrm{PGL}}                  \newcommand\PGammaL{\mathrm{P\Gamma L}}
                  
\newcommand\PSL{\mathrm{PSL}}

                    \newcommand\Sp{\mathrm{Sp}}

\newcommand\AS{\mathrm{AS}}  \newcommand\TW{\mathrm{TW}}    \newcommand\PA{\mathrm{PA}}
\newcommand\HA{\mathrm{HA}}

\newtheorem{theorem}{Theorem}[section]%
\newtheorem{lemma}[theorem]{Lemma}%
\newtheorem{proposition}[theorem]{Proposition}%
\newtheorem{question}[theorem]{Question}%
\newtheorem{hypothesis}[theorem]{Hypothesis}%

\begin{document}

\title[Finite $s$-geodesic transitive graphs under certain girths]
{Finite $s$-geodesic transitive graphs under certain girths}

\thanks{2020 MR Subject Classification 20B15, 20B30, 05C25.}

\author[J.-J. Huang]{Jun-Jie Huang}
\address{Jun-Jie Huang\\
School of Mathematical Sciences, Laboratory of Mathematics
and Complex Systems, MOE\\
Beijing Normal University\\
Beijing \\
100044, P. R. China}
{\email{jjhuang@bnu.edu.cn(J.-J. Huang)}}
\maketitle

\begin{abstract}
  For an integer $s\geq1$ and a graph $\Gamma$, a path $(u_0, u_1, \ldots, u_{s})$ of vertices of $\Gamma$ is called an {\em $s$-geodesic} if it is a shortest path from $u_0$ to $u_{s}$. We say that $\Gamma$ is {\em $s$-geodesic transitive} if, for each $i\leq s$, $\Gamma$ has at least one $i$-geodesic, and its automorphism group is transitive on the set of $i$-geodesics. In 2021, Jin and Praeger [J. Combin. Theory Ser. A 178 (2021) 105349] have studied $3$-geodesic transitive graphs of girth $5$ or $6$, and they also proposed to the problem that to classify $s$-geodesic transitive graphs of girth $2s-1$ or $2s-2$ for $s=4, 5, 6, 7, 8$. The case of $s = 4$ was investigated in [J. Algebra Combin. 60 (2024) 949--963]. In this paper, we study such graphs with $s\geq5$. More precisely, it is shown that a connected $(G,s)$-geodesic transitive graph $\Gamma$ with a nontrivial intransitive normal subgroup $N$ of $G$ which has at least $3$ orbits, where $G$ is an automorphism group of $\Gamma$ and $s\geq 5$, either $\Gamma$ is the Foster graph and $\Gamma_N$ is the Tutte's $8$-cage, or $\Gamma$ and $\Gamma_N$ have the same girth and $\Gamma_N$ is $(G/N,s)$-geodesic transitive. Moreover, it is proved that if $G$ acts quasiprimitively on its vertex set, then $G$ is an almost simple group, and if $G$ acts biquasiprimitively, the stabilizer of biparts of $\Gamma$ in $G$ is an almost simple quasiprimitive group on each of biparts. In addition, $G$ cannot be primitive and biprimitive.
\end{abstract}

\qquad {\textsc k}{\scriptsize \textsc {eywords.}  $s$-geodesic transitive graph, $s$-arc transitive graph, quasiprimtive group, automorphism group, girth} {\footnotesize}

\section{Introduction}

All graphs considered in this paper are finite, connected, undirected and simple graphs.
For a graph $\Ga$, we use $V(\Ga)$ and $\Aut(\Ga)$ to denote its vertex set and the full automorphism group, respectively.
For a positive integer $s$,
an {\em $s$-arc} of $\Ga$ is a sequence of vertices $(u_0,u_1,\ldots,u_s)$ in $\Ga$ such that $u_i$ is adjacent to $u_{i+1}$ for $0\leq i\leq s-1$ and $u_{j-1}\neq u_{j+1}$ for $1\leq j\leq s-1$.
In particular, $1$-arc is simply called {\em arc}.
Let $G$ be a subgroup of $\Aut(\Ga)$.
Then $\Ga$ is called {\em $(G,s)$-arc transitive} if $\Ga$ has at least one $s$-arc and $G$ is transitive on the set of $s$-arcs of $\Ga$,
and a $(G,s)$-arc transitive graph is said to be {\em $s$-arc transitive} if $G=\Aut(\Ga)$.
Moreover, we say that $\Ga$ is {\em $s$-transitive} if $\Ga$ is $s$-arc transitive but not $(s+1)$-arc transitive.
The investigation of $s$-arc transitive graphs was initiated by Tutte~\cite{Tutte1947}, who proved that there are no finite $6$-arc transitive cubic graphs.
More than twenty years later, Weiss~\cite{Weiss81} further showed that there are no finite $8$-arc transitive graphs with valency at least three.
Henceforth, the study of $s$-arc transitive graphs has attracted great attention in algebra graph theory,
see~\cite{Li2001,Li01,LP,Primoz,Praeger93,Praeger93-1,Weiss:s-trans,Zhou2023} and the references therein.

Recently, Praeger et al. extended the definition of $s$-arc transitive graphs to $s$-geodesic transitive graphs.
To explain this, we introduce some notations and definitions, as follows.
For a graph $\Ga$, the {\em girth} of $\Ga$ is the length of the shortest cycle in $\Ga$, denoted by $\g_\Ga$.
The {\em distance} of two distinct vertices $u$ and $v$ in $\Ga$ is the length of the shortest path from $u$ to $v$, denoted by $d_\Ga(u,v)$,
and the {\em diameter} $\diam(\Ga)$ of $\Ga$ is the maximum distance between $u$ and $v$ for all $u,v\in V(\Ga)$.
An $s$-arc $(u_0,u_1,\ldots,u_s)$ of $\Ga$ is called an {\em $s$-geodesic} if $d_\Ga(u_0,u_s)=s$.
Then $\Ga$ is said to be {\em $(G,s)$-geodesic transitive} if for each $1\leq i\leq s\leq\diam(\Ga)$,
$G$ is transitive on the set of all $i$-geodesics of $\Ga$.
Moreover, if $\Ga$ is $(\Aut(\Ga),s)$-geodesic transitive, then it is referred to as {\em $s$-geodesic transitive}.
When $s=\diam(\Ga)$, an $s$-geodesic transitive graph is simply called {\em geodesic transitive}.
Obviously, every $s$-geodesic is also an $s$-arc.
However, some $s$-arcs may not be $s$-geodesics.
For example, all $s$-arcs that lie in a cycle of length $2s-1$ are not $s$-geodesics.
Consequently, the class of $s$-arc transitive graphs is properly included in the class of $s$-geodesic transitive graphs.
It is noted that there are $s$-geodesic transitive graphs that are not $s$-arc transitive, see~\cite{Huang,HFZY,HFZ2025,JDLP} for example.
In addition, the difference from the $s$-arc transitivity is that there is no upper bound on $s$ for $s$-geodesic transitivity, refer to~\cite[Theorem 1.1]{JDLP}.
In view of the above, when studying $s$-geodesic transitive graphs, the most interesting aspect is to study those graphs that are not $s$-arc transitive.
In this paper, we will study such $s$-geodesic transitive graphs.

As is well-known, for an $s$-arc transitive graph $\Ga$, the relation between $s$ and $\g_\Ga$ is $\g_\Ga\geq 2s-2$ (see~\cite[Proposition 17.2]{Biggs}).
Based on this, Jin and Praeger~\cite[Remark 1.1]{JP} suggest that when studying $s$-geodesic transitive but not $s$-arc transitive graphs,
one only needs to study $s$-geodesic transitive graphs with girth $2s-1$ or $2s-2$ and $s\leq 8$.
Then they posed the following question (see~\cite[Problem 1.2]{JP}).

\begin{question}\label{question1}
For $s\in\{4, 5, 6, 7, 8\}$, classify the finite $(G,s)$-geodesic-transitive graphs of girth $2s-1$ or $2s-2$, which are not $(G,s)$-arc-transitive, where $G\leq\Aut(\Ga)$.
\end{question}

We remark that such graphs for $s\leq 4$ have been extensively studied over the years.
For $s=2$, Huang et al.~\cite{HFZY} classified $2$-geodesic transitive graphs of order $p^n$ with $p$ a prime and $n\leq 3$.
Moreover, the $2$-geodesic transitive graphs of girth $3$ and prime valency were completely classified in~\cite[Theorem 1.3]{DJLP15}.
For more results regarding $2$-geodesic transitive graphs of girth $3$,
we refer the readers to~\cite{DJLP13-2,HFZY2025,Jin2015-2,JDLP} and the references therein.
For $s=3$, such graphs of valency at most $5$ are classified in~\cite{HFZ2025,Jin15,Jin18}.
Recently, by using the normal quotient strategy,
Jin, Praeger and Tan~\cite{JP,JT24} presented an interesting and fascinating reduction theorem
for studying $s$-geodesic transitive graphs of girth $2s-1$ or $2s-2$ with $s=3$ or $4$.
More precisely, they showed that for such a graph $\Ga$,
if the girth of the quotient graph $\Ga_N$ of $\Ga$ is not $\g_\Ga$,
where $N$ is a normal subgroup of $G$ with at least $3$ orbits on its vertex set,
then either $(\Ga,\Ga_N)$ are determined in~\cite[Theorem 1.5]{JP} and~\cite[Theorem 1.3]{JT24},
or $\Ga_N$ is a geodesic transitive graph of diameter $3$ and of girth $5$ or $6$.

Building on this work, we give a reduction result for the family of $s$-geodesic transitive graphs of girth $2s-1$ or $2s-2$ with $s\geq5$.
Additionally, we describe various possibilities for the girth and diameter of the quotient graphs.

\begin{theorem}\label{Thm-1}
Let $\Ga$ be a connected $(G,s)$-geodesic transitive graph with girth $2s-2$ or $2s-1$, where $G\leq\Aut(\Ga)$ and $5\leq s\leq 8$.
Let $N$ be a nontrivial normal intransitive subgroup of $G$ that has at least $3$ orbits on $V(\Ga)$.
Then $\Ga$ is an $N$-normal cover of $\Ga_N$ and $s\neq7$.
Moreover, either
\begin{enumerate}[\rm (1)]
  \item $\g_{\Ga_N}=\g_\Ga$, $\diam(\Ga_N)> s$ and $\Ga_N$ is $(G/N,s)$-geodesic transitive; or
  \item $\Ga$ is the Foster graph and the triple $(s,\Ga_N,\g_\Ga)=(6,\Delta_{4,2},10)$.
\end{enumerate}
\end{theorem}

\medskip
\noindent{\bf Remark:} (1) The graph $\Delta_{4,2}$ is the generalized $3$-gon of order $30$ (see Section~\ref{Prel}), and it is also known as the Tutte's $8$-cage.

(2) For a connected $(G,s)$-geodesic transitive graph $\Ga$ with girth $2s-2$ or $2s-1$, where $G\leq\Aut(\Ga)$ and $s\geq2$,
suppose that $\Ga$ is a cover of $\Ga_N$ for some nontrivial intransitive normal subgroup $N$ of $G$, and that $\g_\Ga=\g_{\Ga_N}$.
When $s\leq 4$, there are numerous examples of such $\Ga$ and $\Ga_N$.
For example, the Cayley graph $\Cay(G,S)$ is a cover of $\Cay(G/N,S/N)\cong\K_{p^2}$,
where $G=\la a,b,c\mid a^p=b^p=c^p=1,[a,b]=c,[a,c]=[b,c]=1\ra$, $S=\la a^i,b^i\mid 1\leq i\leq p-1\ra$ and $N=\la c\ra$ with $p$ an odd prime.
These two graphs have girth $3$ and $\Cay(G,S)$ is $2$-geodesic transitive, see~\cite[Example 3.4]{HFZY}.
When $s=3$, the only known examples are $(\Ga,\Ga_N)=(\text{Tetra2AT[165;1],Tetra2AT[55;1]})$, or (Tetra2AT[600;3],Tetra2AT[300;1]).
Such graphs are $3$-geodesic transitive of girth $5$, see~\cite{HFZ2025,Primoz}.
Moreover, examples for $s=4$, we refer to~\cite{ILPP} and~\cite[Theorem 1.3(2)]{JT24}.
When $s\geq5$, we are unaware of any examples of such graphs $\Ga$ and $\Ga_N$.

(2) The are many examples of $s$-geodesic transitive graphs with girth $2s-1$ or $2s-2$ when $s\leq 4$,
see~\cite{DJLP13-2,DJLP15,HFZY,HFZ2025,Jin15,Jin2015-2,Jin18,JDLP,JP,JT24} for example.
However, the only $s$-geodesic transitive graphs with girth $2s-2$ or $2s-1$ and $s\geq5$ that have been discovered so far is the Foster graph.
\medskip

According to Theorem~\ref{Thm-1}, when investigating $(G,s)$-geodesic transitive graphs of girth $2s-2$ or $2s-1$ with $s\geq5$ and $G\leq\Aut(\Ga)$,
the key lies in determining such graphs when $G$ is quasiprimitive or biquasiprimitive on its vertex set.
Our second theorem can reduce the graphs with the aforementioned properties to a single case, and shows that $G$ cannot be primitive or biprimitive on its vertex set.
Before proceeding, we introduce the following definitions.

Let $X$ be a transitive permutation group $X$ on a set $\Omega$.
Then $X$ is called {\em primitive} if it has only trivial blocks in $\Omega$
(a {\em block} is a non-empty subset $\Delta$ of $\Omega$ such that $\Delta^g\cap \Delta\not=\emptyset$ implies $\Delta^g=\Delta$ for every $g\in G$
and a block $\Delta$ is {\em trivial} if $|\Delta|=1$ or $|\Delta|=|\Omega|$),
and {\em quasiprimitive} if every non-trivial normal subgroup of $G$ is transitive on $\Omega$.
Moreover, $X$ is called {\em biquasiprimitive} if each of its nontrivial normal subgroups has at most two orbits and at least one has exactly two orbit,
and called {\em biprimitive} if $G$ has an invariant partition $\{\Omega_1, \Omega_2\}$ such that $G_{\Omega_1}=G_{\Omega_2}$ is primitive on $\Omega_1$ and $\Omega_2$.
Clearly, (bi)primitive group is (bi)quasiprimitive, but the converse is not true.

Let $\Ga$ be a connected $(G,s)$-geodesic transitive graph with valency at least $3$,  where $G\leq\Aut(\Ga)$ and $s\geq5$.
Assume that $\Ga$ has girth $2s-2$ or $2s-1$.
Then $\Ga$ is $(G,4)$-arc transitive by Lemma~\ref{sgeo-to-sarc}.
If $G$ is biquasiprimitive on $V(\Ga)$,
then $\Ga$ is a bipartite graph.
Let $\Delta_1$ and $\Delta_2$ be two biparts of $\Ga$.
Define $G^+=G_{\Delta_1}=G_{\Delta_2}$ and
\begin{equation}\label{eqX}
(X,\Omega)=\left\{
\begin{aligned}
&(G,V(\Ga)),\text{~if~} G \text{~is quasiprimitive on~} V(\Ga),\\
&(G^+,\Delta_1), \text{~if~} G \text{~is biquasiprimitive on~} V(\Ga).
\end{aligned}
\right.
\end{equation}
Next, we will show that $X$ is a quasiprimitive permutation group of type $\AS$ on $\Omega$.

\begin{theorem}\label{Thm:quasi}
Let $\Ga$ be a connected $(G,s)$-geodesic transitive graph of valency at least $3$ and of girth $\g_\Ga=2s-2$ or $2s-1$,
where $G\leq\Aut(\Ga)$ and $s\geq5$.
Suppose that $G$ is quasiprimitive or biquasiprimitive on $V(\Ga)$,
and that $(X,\Omega)$ satisfies Eq.~$\eqref{eqX}$.
Then $X$ is quasiprimitive on $\Omega$ of type $\AS$.
Moreover, $G$ cannot be primitive and biprimitive.
\end{theorem}

We note that we have examples in which $G$ is (bi)quasiprimitive but not (bi)primitive for a graph that meets the conditions of Theorem~\ref{Thm:quasi}.
Moreover, to classify such graph that fulfill Theorem~\ref{Thm:quasi},
it seems difficult to classify almost simple groups with a subgroup satisfies Proposition~\ref{stabilizer}.

The layout of this paper is as follows.
In Section~\ref{Prel}, we list some definitions and results that we need.
The proof of Theorems \ref{Thm-1} and \ref{Thm:quasi} will be presented separately in Sections~\ref{sec:quoitent} and~\ref{sec:quasi}.

\section{Preliminaries}\label{Prel}

In this section, we will list some definitions and properties concerning groups and graphs.
For a positive integer $n$, denote by $\ZZ_n$ the cyclic group of order $n$.
For a prime $p$ and a positive integer $r$, denote by $\ZZ_p^r$ the elementary abelian group of order $p^r$.
For two groups $A$ and $B$, denote by $A\times B$ the direct product of $A$ and $B$, by $A:B$ a semidirect product of $A$ by $B$, and by $A.B$ an extension of $A$ by $B$.
Denote by $\K_{n,n}$ the complete bipartite graph of order $2n$.
For more notations of groups, we refer the readers to Atlas~\cite{Atlas}.

\medskip
Following~\cite{Praeger93,Praeger97}, a finite quasiprimitive permutation group may be divided into eight disjoint types.
In this paper, we provide a simple introduction for quasiprimitive groups of type $\HA$, $\AS$, $\TW$ and $\PA$.
For more details we refer the reader to~\cite{Praeger93,Praeger97}.
Now let $G$ be a finite quasiprimitive permutation group on $\Omega$,
and let $N=\soc(G)$ be the {\em socle} of $G$, that is the product of all minimal normal subgroups of $G$.

\medskip
\noindent
$\mathbf{HA}$ (Holomorph Affine): $N\cong \ZZ_p^k$ is regular on $\Omega$, and $G= N:G_\a$, where $p$ is a prime, $k\geq1$ and $\a\in\Omega$.

\noindent
$\mathbf{AS}$ (Almost Simple): $N=T$ is a non-abelian simple group and $T\unlhd G\leq\Aut(T)$.

\noindent
$\mathbf{TW}$ (Twisted Wreath product): $N=T^k$ acts regularly on $\Omega$, where $T$ is a non-abelian simple group and $k\geq2$.

\noindent
$\mathbf{PA}$ (Product Action): $N=T^k$ is a minimal normal subgroup of $G$ with $k\geq2$ and $T$ a non-abelian simple group.
For $\a\in\Omega$, there exists a non-trivial proper subgroup $R$ of $T$ such that $N_\a$ is a subgroup of $R^k$ which projects surjectively onto each of the
$k$ direct factors $R$.
\medskip

The following well-known result, as detailed in \cite{Weiss81,Weiss:s-trans},
gives the structure of the stabilizer of the automorphism group of $4$-arc transitive graphs.

\begin{proposition}\label{stabilizer}
Let $\Ga$ be a connected $(G,s)$-arc transitive graph, where $s\geq 4$ and $G\leq\Aut(\Ga)$.
Let $q=p^f$ be a prime power and let $\{u,v\}$ be an edge of $\Ga$.
Then the kernel of $G_{uv}$ acting on $\Ga(v)$ is a $p$-group, and one of the following holds:
\begin{enumerate}[\em (i)]
  \item $s=4$ and $G_u\cong [q^2]:(\ZZ_{(q-1)/(3,q-1)}.\PGL(2,q)).o$ with $|o|\bigm|(3,q-1)f$;
  \item $s=5$, $p=2$ and $G_u\cong [q^3]:\GL(2,q).\ZZ_e$ with $e\mid f$;
  \item $s=7$, $p=3$ and $G_u\cong [q^5]:\GL(2,q).\ZZ_e$ with $e\mid f$.
\end{enumerate}
\end{proposition}

Let $\Ga$ be an $s$-geodesic transitive graph with $s\leq\diam(\Ga)$.
Let $u\in V(\Ga)$. For each $1\leq i\leq s$,
every vertex in $\Ga_i(u)$ is adjacent to the same number of other vertices in $\Ga_i(u)$, say $a_i$.
Furthermore, every vertex in $\Ga_i(u)$ is adjacent to the same number of vertices in $\Ga_{i-1}(u)$, say $c_i$.
Similarly, if $i<\diam(\Ga)$, then every vertex in $\Ga_i(u)$ is adjacent to the same number of vertices in $\Ga_{i+1}(u)$, say $b_i$.
If $i=\diam(\Ga)$, then let $b_i=0$.
Clearly, $a_i+b_i+c_i$ is equal to the valency of $\Ga$.
If $\Ga$ is geodesic transitive with diameter $d$,
then the constants are well-defined for all $1\leq i\leq d$,
and the array $\{b_0,b_1,\ldots,b_{d-1};c_1,c_2,\ldots,c_{d}\}$,
where $b_0$ is the valency of $\Ga$, is called the {\em intersection array} of $\Ga$.
These parameters play an important role in the study of $s$-geodesic transitive graphs,
see~\cite{BCN,HFZY,HFZ2025,Jin15,JP,JT24}, for example.

The following conclusion gives some properties of $s$-geodesic transitive graphs, see~\cite[Lemma 2.4]{HFZ2025} and~\cite[Lemma 2.3]{JP}.

\begin{proposition}\label{dis-tran}
Let $s$ be a positive integer, and let $\Ga$ be an $s$-geodesic transitive graph with a vertex $u$.
\begin{enumerate}[\rm (1)]
  \item $|\Ga_i(u)|b_i=|\Ga_{i+1}(u)|c_{i+1}$ for all $1\leq i\leq s-1$.
  \item If $b_s\leq 1$, then $\Ga$ is geodesic transitive.
\end{enumerate}
\end{proposition}

According to~\cite[P. 84]{GR}, a {\em generalized polygon}, more precisely, a {\em generalized $d$-gon}, is a bipartite graph with diameter $d$ and girth $2d$.
Based on the work in~\cite{Weiss85}, we refer to the generalized polygons associated with the simple groups $\PSL_3(q)$, $\Sp_4(q)$ and $\G_2(q)$,
where $q$ is a prime power, as {\em classical generalized polygons}.
These are respectively denoted as $\Delta_{3,q}$, $\Delta_{4,q}$ and $\Delta_{6,q}$.
The following result is quoted from~\cite{BCN,JDLP}
(the first four columns in Table~\ref{tab:geod-sarc} are from~\cite[Table 1]{JDLP}, and the last three columns are from~\cite[Section 6.5]{BCN}).

\begin{proposition}\label{geodesic-arc}
Let $\Ga$ be a connected geodesic transitive graph with valency $k$ and let $s\geq4$.
Then $\Ga$ is $s$-transitive if and only if $\Ga$ is isomorphic to one of the graphs listed in Table~$\ref{tab:geod-sarc}$
$($In Table~$\ref{tab:geod-sarc}$, $d=\diam(\Ga))$.
\end{proposition}

\begin{table}[!hbtp]
\centering
\caption{Geodesic transitive and $s$-transitive graphs with $s\geq4$} \label{tab:geod-sarc}
\begin{tabular}{l|l|l|l|l|l|l} \hline
  $\Ga$  & $k$ &  $s$ &$d$ & $\g_\Ga$ & Intersection array & $\Aut(\Ga)$ \\ \hline
  the Foster graph   & $3$ & $5$  & $8$ & $10$ & $\{3,2,2,2,2,1,1,1;$ & $\ZZ_3.\PGammaL(2,9)$\\
                                       &     &       &     &      & $~~~~~~1,1,1,1,2,2,2,3\}$ & \\
  the Biggs-Smith graph & $3$ & $4$ & $7$ & $9$  & $    \{3,2,2,2,1,1,1;$  & $\PSL(2,17)$\\
                                        &     &       &     &     & $~~~~~~1,1,1,1,1,1,3\}$  & \\
  $\Delta_{3,q}$, $q$ a prime power & $q+1$ & $4$  & $3$ &  $6$ & $\{q+1,q,q;1,1,q+1\}$ & $\Aut(\PSL(3,q))$ \\
  $\Delta_{4,q}$, $q=2^m$ and $m\geq1$ & $q+1$ & $5$  & $4$ & $8$ & $\{q+1,q,q,q;$  & $\Aut(\Sp(4,q))$\\
                                         &     &       &     &     & $~~~~~~1,1,1,q+1\}$  &  \\
  $\Delta_{5,q}$, $q=3^m$ and $m\geq1$ & $q+1$ & $7$  & $6$ & $12$&  $\{q+1,q,q,q,q,q;$ & $\Aut(\G_2(q))$ \\
                                         &     &       &     &     & $~~~~~~1,1,1,1,1,q+1\}$  &  \\\hline
\end{tabular}
\end{table}

Let $\Ga$ be a graph with $G\leq \Aut(\Ga)$, and let $N$ be a normal subgroup of $G$ such that $N$ is intransitive on $V(\Ga)$.
The {\it normal quotient graph} $\Ga_N$ of $\Ga$ induced by $N$ is defined to be the graph
with vertex set the orbits of $N$ and two orbits are adjacent in $\Ga_N$
if $\Ga$ has an edge incident to some vertices in the two orbits, respectively.
Furthermore, $\Ga$ is called a {\it normal $N$-cover} or simply {\it cover} of $\Ga_N$
if for every vertex $v\in V(\Ga)$ and the orbit $O(v)$ of $N$ containing $v$,
$v$ and $O(v)$ has the same valency in $\Ga$ and $\Ga_N$, respectively.
The following proposition gives a reduction for studying $(G,s)$-geodesic transitive graphs, refer to~\cite[Lemma 3.2]{JP}.

\begin{proposition}\label{redu:geodes}
Let $\Ga$ be a connected $(G,s)$-geodesic transitive graph, where $s\geq2$ and $G\leq\Aut(\Ga)$.
Let $1\neq N\lhd G$ be intransitive on $V(\Ga)$.
Suppose that $\Ga$ is not a complete multipartite graph.
Then either
\begin{enumerate}[\rm (i)]
  \item $N$ has $2$ orbits on $V(\Ga)$ and $\Ga$ is bipartite; or
  \item $N$ has at least $3$ orbits on $V(\Ga)$, $N$ is semiregular on $V(\Ga)$, $\Ga$ is a cover of $\Ga_N$ and $\Ga$ is $(G/N,s')$-geodesic transitive where $s'=\min\{s,\diam(\Ga_N)\}$.
\end{enumerate}
\end{proposition}

\section{Proof of Theorem~\ref{Thm-1}}\label{sec:quoitent}

In this section, we will prove Theorem~\ref{Thm-1} through a series of lemmas.
Firstly, we present the following hypothesis.

\begin{hypothesis}\label{hypothesis}
Let $\Ga$ be a connected $(G,s)$-geodesic transitive graph with girth $\g_\Ga=2s-2$ or $2s-1$, where $G\leq\Aut(\Ga)$  and $s\geq2$.
Let $1\neq N\lhd G$ be an intransitive subgroup with at least $3$ orbits on $V(\Ga)$ and let $\Sig=\Ga_N$.
Suppose that $\Ga$ is a cover of $\Sig$, and that $k=\g_\Sig$.
\end{hypothesis}

The following two lemmas aim to study the bounds of the girth of $\Sig$ and the distances between some vertices in $\Ga$ or $\Sig$.

\begin{lemma}\label{proper:girth}
Suppose Hypothesis~$\ref{hypothesis}$.
Then $2s-4 \leq k\leq \g_\Ga$, and $\Sig$ is $(s-1)$-transitive with $s\neq 7$.
\end{lemma}

\proof Notice that $\Ga$ is $(G,s-1)$-geodesic transitive graph.
Since $\g_\Ga=2s-1$ or $2s-2$, all $(s-1)$-arcs of $\Ga$ are $(s-1)$-geodesics,
and hence $\Ga$ is also $(G,s-1)$-transitive.
This property, together with $\Ga$ is a cover of $\Sig$, implies that $\Sig$ is $(G/N,s-1)$-transitive (see~\cite[Lemma 2.5]{LP}).
According to~\cite[Proposition 17.2]{Biggs}, we have $k\geq 2s-4$.
Moreover, since $\Ga$ is a cover of $\Sig$, each $\g_\Ga$-cycle of $\Ga$ induces a $\g_\Ga$-cycle of $\Sig$, and so $k\leq \g_\Ga$.
Therefore, $2s-4 \leq k\leq \g_\Ga$.

If $s=7$, then $\Sig$ is $6$-arc transitive, so is $7$-arc transitive (see~\cite[Theorem]{Weiss81}).
Within $\Sig$, there are two distinct types of $7$-arcs.
One type is the $7$-geodesic, and the other type is the $7$-arcs that lies in some $12$-cycles or $13$-cycles.
This co-existence of two different types of $7$-arcs contradicts the $7$-arc transitivity of $\Sig$.
Thus, $s\neq7$, which concludes the proof. \qed

\begin{lemma}\label{proper:distance}
Suppose Hypothesis~$\ref{hypothesis}$.
Assume that $(B_1,B_2,\cdots,B_k)$ is a $k$-cycle of $\Sig$ with $k<\g_\Ga$.
Then the following holds.
\begin{enumerate}[\rm (1)]
  \item $d_\Sig(B_1,B_{i})=i-1$ for all $2\leq i\leq r+1$ and $d_\Sig(B_1,B_{r+2})=r-1$ or $r$ depending on $k$ is even or odd, respectively, where $r=[k/2]$.
  \item There is a $k$-arc $(u_1,u_2,\ldots,u_k,u_1')$ of $\Ga$, where $u_i\in B_i$ with $1\leq i\leq k$ and $u_1'\in B_1\setminus\{u_1\}$ such that
   $d_\Ga(u_1,u_i)=i-1$ for all $2\leq i\leq s$.
\end{enumerate}
\end{lemma}

\proof Let $r=[k/2]$.
Since $\Sig$ has girth $k$, we have $d_\Sig(B_1,B_{i})=i-1$ for all $2\leq i\leq r+1$.
As part (1) states, when $i = r + 2$, $d_{\Sig}(B_1,B_{r+2})$ is equal to $r - 1$ if $k$ is even, and $r$ if $k$ is odd.

Since $\Ga$ is a cover of $\Sig$,
there exist vertices $u_i\in B_i$ for $1\leq i\leq k$ and $u_1'\in B_1$ such that $(u_1,u_2,\ldots,u_k,u_1')$ forms a $k$-arc of $\Ga$.
In particular, $u_1\neq u_1'$ as $\g_\Ga>k$.
For each $2\leq i\leq s$, we have $d_\Ga(u_1,u_i)=\ell_i\leq i-1$.
Assume that $\ell_i<i-1$.
Then there exists an $\ell_i$-geodesic of $\Ga$ passing through $u_1$ and $u_i$.
This $\ell_i$-geodesic, combined with the $i$-arc $(u_1,\ldots,u_i)$, forms a loop of length $\ell_i+i-1<2i-2\leq 2s-2\leq \g_\Ga$,
which is clearly impossible.
Therefore, $\ell_i=i-1$, which completes the proof. \qed

Next, we will determine the quotient graph $\Sig$, except for the case where $\g_\Sig=\g_\Ga$.

\begin{lemma}\label{reduce}
Suppose Hypothesis~$\ref{hypothesis}$.
If $5\leq s\leq 8$, then $s\neq7$, and either
\begin{enumerate}[\rm (1)]
  \item $\g_{\Sig}=\g_\Ga$, $\diam(\Sig)> s$ and $\Sig$ is $(G/N,s)$-geodesic transitive; or
  \item $\Sig$ is geodesic transitive, and $(s,\Sig,\g_\Ga)=(5,\Delta_{3,q},8)$, $(6,\Delta_{4,q},10)$ or $(8,\Delta_{5,q},14)$.
\end{enumerate}
\end{lemma}

\proof Since $\Ga$ is a cover of $\Sig$,
by Proposition~\ref{redu:geodes},
$\Sig$ is $(G/N,s')$-geodesic transitive with $s'=\min\{s,\diam(\Sig)\}$.
By Lemma~\ref{proper:girth}, $\Sig$ is $(s-1)$-transitive with $s\neq 7$.

Suppose that $\diam(\Sig)\leq s$.
Then $s'=\diam(\Sig)$, and hence $\Sig$ is geodesic transitive.
Note that the Foster graph is a $5$-transitive graph with diameter $8$,
and the Biggs-Smith graph is a $4$-transitive graph with diameter $7$ (see Table~\ref{tab:geod-sarc}).
Since $s\geq5$ and $\diam(\Sig)\leq s$,
it follows from Proposition~\ref{geodesic-arc} that $(s,\Sig)=(5,\Delta_{3,q})$, $(6,\Delta_{4,q})$ or $(8,\Delta_{5,q})$.
In this case, $\Sig$ is a bipartite graph,
and so $\Ga$ is bipartite, yielding $\g_\Ga=2s-2$.
Thus, part (2) holds.

Suppose that $\diam(\Sig)> s$.
Then $\Sig$ is $s$-geodesic transitive.
Recall that $k=\g_\Sig$.
According to Lemma~\ref{proper:girth}, we have the inequality $2s-4\leq k\leq \g_\Ga$.
To finish the proof, suppose that $k<\g_\Ga$.
Then either $k=2s-4$ or $2s-3$, or $k=2s-2$ and $\g_\Ga=2s-1$.
Let $(B_1,\ldots,B_k)$ be a $k$-cycle of $\Sig$.
By Lemma~\ref{proper:distance}, there exist vertices $u_i\in B_i$ such that $(u_1,\ldots,u_k,u_1')$ is a $k$-arc,
and $d_\Ga(u_1,u_j)=j-1$ for all $2\leq j\leq s$, where $u_1'\in B_1$.

Let $k=2s-4$ or $2s-3$.
Then $(B_1,B_2,\ldots,B_s)$ is an $(s-1)$-arc of $\Sig$ with $d_\Sig(B_1,B_s)<s-1$.
Since $\diam(\Sig)> s$, $\Sig$ has an $(s-1)$-geodesic, say $(B_1,B_2',B_3',\ldots,$$ B_s')$.
However, no element of $\Aut(\Sig)$ can map the $(s-1)$-arc $(B_1,B_2,\ldots,B_s)$ to the $(s-1)$-geodesic $(B_1,B_2',B_3',\ldots, B_s')$,
contradicting to $\Sig$ is $(s-1)$-arc transitive.

Let $k=2s-2$ and $\g_\Ga=2s-1$.
Then $d_\Sig(B_1,B_j)=d_\Ga(u_1,u_j)=j-1$ for every $2\leq j\leq s$, and $d_\Sig(B_1,B_{s+1})=s-2$.
Thus, $d_\Ga(u_1,u_{s+1})\leq s$.
Assume that $d_\Ga(u_1,u_{s+1})=\ell$ with $\ell\leq s-1$.
Since $\Ga$ is $(G,s)$-geodesic transitive, $G_{u_1}$ is transitive on $\Ga_\ell(u_1)$,
and so $u_{s+1}^x=u_{\ell+1}$ for some $x\in G_{u_1}$.
It follows that $B_1^x=B_1$ and $B_{s+1}^x=B_{\ell+1}\in\Sig_\ell(B_1)$.
In particular, $\ell=s-2$.
Let $(u_1,w_2,\ldots,w_{s-2},u_{s+1})$ be an $(s-2)$-geodesic of $\Ga$.
Then the sequence $(u_1,w_2,\ldots,w_{s-2},u_{s+1},u_s,\ldots,u_2)$ forms a loop of length $2s-2$,
which contradicts the fact that $\g_\Ga=2s-1$.
Therefore, $d_\Ga(u_1,u_{s+1})= s$.

Recall that $\diam(\Sig)> s$.
Then $\Sig$ has an $s$-geodesic, say $(B_1,C_2,\ldots,C_{s+1})$.
Since $\Ga$ is a cover of $\Sig$, there exist vertices $v_j\in C_j$ ($2\leq j\leq s+1$) such that $(u_1,v_2,\ldots,v_{s+1})$ is an $s$-geodesic of $\Ga$.
By the $(G,s)$-geodesic transitivity of $\Ga$, there is an element that maps $(u_1,u_{s+1})$ to $(u_1,v_{s+1})$,
and thus it also maps $(B_1,B_{s+1})$ to $(B_1,C_{s+1})$.
However, $d_\Sig(B_1,B_{s+1})=s-2$ and $d_\Sig(B_1,C_{s+1})=s$, which is impossible.
This completes the proof. \qed

To avoid confusion, for a graph $\Ga$,
we use $a_i^\Ga$, $b_i^\Ga$ and $c_i^\Ga$ to denote the parameters $a_i$, $b_i$ and $c_i$ respectively, which are defined in Section~\ref{Prel}.
In the next lemma, we will determine the graph $\Ga$ when $\Ga$ is a cover of $\Sig$ and $\Sig$ satisfies case (2) of Lemma~\ref{reduce}.

\begin{lemma}\label{geodesic-cover}
Let $\Ga$ be a connected $(G,s)$-geodesic transitive graph, where $G\leq \Aut(\Ga)$.
Then $\Ga$ is a cover of a graph $\Sig$ with the triples $(s,\Sig,\g_\Ga)=(5,\Delta_{3,q},8)$, $(6,\Delta_{4,q},10)$ or $(8,\Delta_{5,q},14)$
if and only if $\Ga$ is the Foster graph and $(s,\Sig,\g_\Ga)=(6,\Delta_{4,2},10)$.
\end{lemma}

\proof Notice that the Foster graph is a geodesic transitive graph of diameter $8$ and of girth $10$ (see Table~\ref{tab:geod-sarc}).
From~\cite[P. 398]{BCN} we know that the Foster graph is the unique $3$-cover of $\Delta_{4,2}$.
Therefore, the sufficiency part holds.

To prove the necessity,
assume that $\Ga$ is a cover of $\Sig$ with the triples $(s,\Sig,\g_\Ga)=(5,\Delta_{3,q},8)$, $(6,\Delta_{4,q},10)$ or $(8,\Delta_{5,q},14)$.
By Table~\ref{tab:geod-sarc}, $\Sig$ is a bipartite graph satisfying
\begin{equation}\label{eqaiSig}
\diam(\Sig)=s-2,~\g_\Sig=2s-4 \text{~and~} c_i^\Sig=1 \text{~for all~} 1\leq i<\diam(\Sig).
\end{equation}
Thus, $\Ga$ is a bipartite graph of girth $2s-2$.
Let $(u_0,u_1,\ldots,u_{2s-3})$ be a cycle of length $2s-2$.
Then $d_\Ga(u_0,u_i)=i$ for each $1\leq i\leq s-1$.
Since $\Ga$ is a cover of $\Sig$,
we may let $B_i\in V(\Sig)$ be such that $u_i\in B_i$, where $0\leq i\leq 2s-3$.
Then $d_\Sig(B_0,B_j)=j$ for $1\leq j\leq s-2$,
and $d_\Sig(B_0,B_{s-1})=s-3$ or $s-2$ as $\diam(\Sig)=s-2$ (also see Lemma~\ref{proper:distance}).
If $d_\Sig(B_0,B_{s-1})=s-2$, then $\Sig$ must contain a cycle of length odd,
which contradicts the fact that $\Sig$ is a bipartite graph.
Thus, $d_\Sig(B_0,B_{s-1})=s-3$.

Since $\g_\Sig=2s-4$, there exists a cycle of length $2s-4$ passing through the $s$-arc $(B_0,\ldots,B_{s-1})$ of $\Sig$.
Let $C_j\in V(\Sig)$, where $1\leq j\leq s-4$, be such that
\begin{equation*}
(B_0,\ldots,B_{s-1},C_{s-4},C_{s-5},\ldots,C_1)
\end{equation*}
is a $(2s-4)$-cycle of $\Sig$.
It follows from $\Ga$ is a cover of $\Sig$ that
\begin{equation*}
(u_0,u_1,\ldots,u_{s-1},w_{s-4},w_{s-5},\ldots,w_1,u_0')
\end{equation*}
is a $(2s-3)$-arc of $\Ga$, where $w_j\in C_j$ and $u_0'\in B_0$.
In particular, $d_\Sig(B_0,C_{s-4})=s-4$ and $u_0'\neq u_0$.
Since $\Ga$ is a bipartite graph and $d_\Ga(u_0,u_{s-1})=s-1$, we have $d_\Ga(u_0,w_{s-4})=s$ or $s-2$.
Note that $\Ga$ is $(G,s)$-geodesic transitive, and then $G_{u_0}$ is transitive on $\Ga_j(u_0)$ for all $1\leq j\leq s$.
If $d_\Ga(u_0,w_{s-4})=s-2$, then $w_{s-4},u_{s-2}\in\Ga_{s-2}(u_0)$, and so $w_{s-4}^x=u_{s-2}$ for some $x\in G_{u_0}$.
Thus, $(B_0,C_{s-4})^x=(B_0,B_{s-1})$, which is impossible because $d_\Sig(B_0,C_{s-4})=s-4$ and $d_\Sig(B_0,B_{s-2})=s-2$.
We therefore conclude that $d_\Ga(u_0,w_{s-4})=s$

Assume that $b_{s-1}^\Ga\geq2$.
Then there exists $w\in\Ga(u_{s-1})\cap\Ga_{s}(u_0)$ such that $w\neq w_{s-4}$.
Since $G_{u_0}$ is transitive on $\Ga_s(u_0)$, we have $w^y=w_{s-4}$ for some $y\in G_{u_0}$.
Let $C\in V(\Sig)$ contain the vertex $w$.
If $C=C_{s-4}$, then $w,w_{s-4}\in\Ga(u_{s-1})\cap C$.
However, $\Ga$ is a cover of $\Sig$ means that $|\Ga(u_{s-1})\cap C|=1$, a contradiction.
Thus, $C\neq C_{s-4}$.
Since $(B_0,C)^y=(B_0,C_{s-4})$, we obtain $C\in \Sig_{s-4}(B_0)$,
and so $\Sig(B_{s-1})\cap\Sig_{s-4}(B_0)$ contains at least two blocks $C$ and $C_{s-4}$.
This implies that $c_{s-4}^\Sig=|\Sig(B_{s-1})\cap\Sig_{s-4}(B_0)|\geq2$, contradicting to Eq.~\eqref{eqaiSig}.

The argument of the previous paragraph shows that $b^\Ga_{s-1}=1$.
From Proposition~\ref{dis-tran} we know that $\Ga$ is geodesic transitive.
Recall that $s\geq5$ and $\Ga$ is $(G,s)$-geodesic transitive of girth $2s-2$.
Then $\Ga$ is also $4$-arc transitive, and so Proposition~\ref{geodesic-arc} is applicable.
By Table~\ref{tab:geod-sarc}, we have either $\Ga=\Sig$, or $\Ga$ is the Foster graph or the Biggs-Smith graph.
Since $2s-2=\g_\Ga\neq\g_\Sig$ and the Biggs-Smith graph has girth $9$, we conclude that $\Ga$ is the Foster graph.
This completes the proof.\qed

We are read to prove Theorem~\ref{Thm-1}.

\medskip
\noindent{\bf Proof of Theorem~\ref{Thm-1}.}
Let $\Ga$ be a $(G,s)$-geodesic transitive graph with girth $\g_\Ga=2s-2$ or $2s-1$, where $G\leq\Aut(\Ga)$  and $s\geq5$.
Let $1\neq N\lhd G$ be an intransitive subgroup that has at least $3$ orbits on $V(\Ga)$ and let $\Sig=\Ga_N$.
Since $\g_\Ga\geq8$, we conclude that $\Ga$ cannot be a compete multipartite graph.
By Proposition~\ref{redu:geodes}, $\Ga$ is a cover of $\Sig$, and so the Hypothesis~\ref{hypothesis} is satisfied.
Combining Lemmas~\ref{reduce} and~\ref{geodesic-cover}, we obtain Theorem~\ref{Thm-1}. \qed

\section{Proof of Theorem~\ref{Thm:quasi}}\label{sec:quasi}

For a connected $(G,s)$-geodesic transitive graph $\Ga$ with $s\geq2$ and $G\leq\Aut(\Ga)$,
we know that all $s$-arcs of $\Ga$ are $s$-geodesics,
and hence the $(G,s)$-geodesic transitivity of $\Ga$ means that $\Ga$ is also $(G,s)$-arc transitive.
This fact will be used repeatedly, we record it here.

\begin{lemma}\label{sgeo-to-sarc}
Let $\Ga$ be a connected $(G,s)$-geodesic transitive graph with $s\geq2$ and $G\leq\Aut(\Ga)$.
If $\Ga$ has girth at least $2s$, then $\Ga$ is $(G,s)$-arc transitive.
\end{lemma}

In the following, we will prove Theorem~\ref{Thm:quasi} with the help of three lemmas.

\begin{lemma}\label{reduceAS}
Let $\Ga$ be a connected $(G,s)$-geodesic transitive graph of valency of girth $\g_\Ga=2s-2$ or $2s-1$,
where $G\leq\Aut(\Ga)$ and $s\geq5$.
Suppose that $G$ is quasiprimitive or biquasiprimitive on $V(\Ga)$,
and that the pair $(X,\Omega)$ satisfies Eq.~$\eqref{eqX}$.
Then $X$ is quasiprimitive on $\Omega$ of type $\AS$.
\end{lemma}

\proof Clearly, $X$ is quasiprimitive on $\Omega$ if $G$ is quasiprimitive on $V(\Ga)$.
Assume that $G$ is biquasiprimitive on $V(\Ga)$.
Then $\Ga$ is a bipartite graph.
Let $\Delta_1$ and $\Delta_2$ be the two bipartite sets of $\Ga$.
Recall that $X=G^+=G_{\Delta_1}=G_{\Delta_2}$.
Then $X$ is a normal subgroup of $G$ with index $2$.

Let $K$ be the kernel of the action of $X$ on $\Delta_1$.
Then $K$ is normal in $X$ but not normal in $G$.
Thus, there is an element $x\in G\setminus X$ such that $K^x\neq K$.
This implies that $G=\la X,x\ra$ and $x^2\in G$.
Let $L=K^x$.
Then $L^x=K^{x^2}=K$ and $\Delta_1^x=\Delta_2$,
and so $L$ is the kernel of $X$ acting on $\Delta_2$.
Let $u\in\Delta_1$.
Then $X_u=G_u$.
If $K\neq 1$, then $L\neq1$, and by the connectivity of $\Ga$,
we obtain $1\neq K^{\Ga(u)}\unlhd X_u^{\Ga(u)}=G_u^{\Ga(u)}$.
Since $\Ga$ is $(G,2)$-arc transitive, $G_u^{\Ga(u)}$ is $2$-transitive,
implying that $K$ is transitive on $\Ga(u)$.
Similarly, $L$ is transitive on $\Ga(v)$ for some $v\in\Ga(u)$.
Therefore, the subgraph $[\Ga(u),\Ga(v)]\cong\K_{n,n}$ with $n=|\Ga(u)|$.
Again by the connectivity of $\Ga$, we have $\Ga\cong\K_{n,n}$.
However, $\K_{n,n}$ is not $4$-arc transitive, which is impossibly by Lemma~\ref{sgeo-to-sarc}.

We therefore conclude that $X$ is faithful on $\Delta_1$ and $\Delta_2$.
Since $G$ is biquasiprimitive on $V(\Ga)$, $X\cong X^{\Delta_1}\cong X^{\Delta_2}$ is quasiprimitive.
Therefore, $X$ is quasiprimitive on $\Omega$.
Next, we will show that $X$ is of type $\AS$.

Let $u\in \Omega$.
Since $\Ga$ is $(G,4)$-arc transitive (see Lemma~\ref{sgeo-to-sarc}),
by~\cite[Theorem 2]{Praeger93} and~\cite[Theorem 2.3]{Praeger93-1},
we conclude that $X$ is of type $\HA$, $\AS$, $\TW$, or $\PA$.
Let $N=\soc(X)$.

Assume that $X$ is of type $\HA$ or $\TW$.
Then $N$ is regular on $\Omega$.
From~\cite[Proposition 2.3]{Li01} and~\cite[Theorem 1.2]{Zhou2023} we know that $\Ga$ is not $4$-arc transitive, which is a contradiction.

Assume that $X$ is of type $\PA$.
Then $N\cong T^k$ and $\Omega=O^k$ for a non-empty set $O$, where $T$ is a nonabelian simple group and $k\geq2$.
In this case, we can set $u=(\a,\a,\ldots,\a)$ with $\a\in O$.
Then $N_u=T_\a^k$ and $1\neq N_u^{\Ga(u)}\unlhd X_u^{\Ga(u)}$.
If $\Ga$ has valency $2$, then $\Ga$ is the cycle of length $2s-1$ or $2s-2$.
In the case, $\Ga$ has no $s$-geodesic, which is clearly impossible.
If $\Ga$ has valency $3$,
then $X_u\cong\Sy_4$ or $\Sy_4\times\ZZ_2$ is solvable (see Proposition~\ref{stabilizer}).
In particular, $X_u^{\Ga(u)}\cong\Sy_3$.
Since $T_\a^k= N_u\unlhd X_u$ with $k\geq2$, we have $T_\a\cong\ZZ_2$ and $k=2$ or $3$,
and so $N_u$ is a $2$-group.
However, $X_u^{\Ga(u)}$ is $2$-transitive means that $N_u^{\Ga(u)}$ is transitive,
and hence $3$ divides $|N_u|$, which is impossible.
Therefore, $\Ga$ has valency at least $4$.
By Proposition~\ref{stabilizer} we know that $\soc(X_u^{\Ga(u)})\cong\PSL(2,q)$,
and $X_u$ has a unique insoluble composition factor, namely, $\PSL(2,q)$.
It then follows that $\PSL(2,q)\leq N_u^{\Ga(u)}$, and $N_u$ has a unique insoluble composition factor.
Since $T_\a^k= N_u\leq X_u$,
$T_\a$ is insoluble, and hence $N_u$ has more than one insoluble composition factor, a contradiction.
Therefore, $X$ is quasiprimitive on $\Omega$ of type $\AS$, completing the proof. \qed

For a $(G,s)$-geodesic transitive graph $\Ga$ with an $s$-geodesic $(u_0,u_1,\ldots,u_{s+1})$,
we know that for each $0\leq i\leq s-1$,
the stabilizer $G_{u_0,u_1,\ldots,u_i}$ is transitive on $\Ga(u_i)\cap\Ga_{i+1}(u_0)$ with degree $b_i$.
Here, $b_0$ is the valency of $\Ga$ and the values of $b_i$ for $1\leq i\leq s-1$ are defined in Section~\ref{Prel}.
It follows that
\begin{equation}\label{eqstab}
b_0b_1\cdots b_{s} \bigm| |G_{u_0}|,
\end{equation}
also see~\cite[Corollary 2.6]{Huang}.
This fact will be repeatedly used in the proof of the following two lemmas.

\begin{lemma}\label{primitive}
Let $\Ga$ be a connected $(G,s)$-geodesic transitive graph with girth $\g_\Ga=2s-2$ or $2s-1$,  where $G\leq\Aut(\Ga)$ and $s\geq5$.
Then $G$ cannot be primitive on $V(\Ga)$.
\end{lemma}

\proof Assume that $G$ is primitive on $V(\Ga)$.
Let $k$ be the valency of $\Ga$ and let $u\in V(\Ga)$.
By Lemma~\ref{sgeo-to-sarc}, $\Ga$ is $(G,s-1)$-transitive, so is $4$-arc transitive.
Since $G$ acts primitively on $V(\Ga)$,
$\Ga$ is completely determined in~\cite[Theorem 1.3]{Li2001}.
This result implies that one of the following holds:
\begin{enumerate}[\rm (i)]
  \item $s=5$, $k=3$, and $G_u\cong\Sy_4\times\ZZ_2$;
  \item $s=6$, and $(k,G_u)=(3,\Sy_4)$, $(5,\ZZ_2^4:\GL(2,4))$, $(6,\ZZ_5^2:\GL(2,5))$ or $(14,\ZZ_{13}^2.(\ZZ_4:\PGL(2,13)))$.
\end{enumerate}
Moreover, each graph listed in~\cite[Theorem 1.3]{Li2001} is not isomorphic to any of the graphs in Table~\ref{tab:geod-sarc}.
In particular, $\Ga$ is not geodesic transitive.
It follows from Proposition~\ref{dis-tran} that $b_{s-1}\neq1$.

Assume that $s=6$. Then $b_5\neq 1$.
Since $\g_\Ga=2s-1$ or $2s-2$, we have
\begin{equation}\label{eq:primitive6}
a_i=0, b_i=k-1 \text{~and~} c_i=1 \text{~for all~} 1\leq i\leq 4; b_5\geq 2, c_5\geq1.
\end{equation}
By Eq.~\eqref{eqstab}, we have $kb_1b_2b_3b_4b_5=48b_5$ divides $|G_u|$, and hence $b_5=1$,
which is impossible.

Assume that $s=5$. Then $b_4\neq 1$.
Again by $\g_\Ga=2s-1$ or $2s-2$, we have
\begin{equation}\label{eq:primitive5}
a_i=0, b_i=k-1 \text{~and~} c_i=1 \text{~for all~} 1\leq i\leq 3; b_4\geq 2, c_4\geq1.
\end{equation}
By Eq.~\eqref{eqstab}, we have that $k(k-1)^3b_4$ is a divisor of $|G_u|$.

Let $(k,G_u)=(3,\Sy_4)$. Then $b_4=1$, a contradiction.

Let $(k,G_u)=(5,\ZZ_2^4:\GL(2,4))$.
Then $b_4\mid 3^2$.
Since $a_4+b_4+c_4=k$ and by Eq.~\eqref{eq:primitive5}, we get  $(a_4,b_4,c_4)=(0,3,1)$.
It follows that $\g_\Ga\geq10$, a contradiction.

Let $(k,G_u)=(6,\ZZ_5^2:\GL(2,5))$.
By~\cite[p. 237 and Remark 5]{SW}, we know that $\g_\Ga=12$ (this can also be verified by Magma~\cite{Magma}), leading to a contradiction.

Let $(k,G_u)=(14,\ZZ_{13}^2.(\ZZ_4:\PGL(2,13)))$.
By Proposition~\ref{dis-tran}, we have $|\Ga_4(u)|=\frac{14\cdot13^3}{c_4}$.
In particular, $c_4$ divides $14\cdot13^3$.
Since $b_4\geq 2$ and $a_4+b_4+c_4=14$, we conclude that $c_4=1,2$ or $7$.
Recall that $\Ga$ is $(G,5)$-geodesic transitive.
Then $G_u$ is transitive on $\Ga_4(u)$.
By~\cite[Table 1]{Li2001}, we find that the group $G$ is the Monster simple group $\mathrm{M}$.
Notice that the permutation representation of $G_u$ is given in the online of Atlas~\cite{AtlasOnline},
and then we can use Magma~\cite{Magma} to obtain some properties of $G_u$.
If $c_4=7$, then $G_u$ has a subgroup of index $|\Ga_4(u)|=2\cdot13^3$.
However, checking by Magma~\cite{Magma}, $G_u$ has no subgroups with the above index, a contradiction.
Thus, $c_4=1$ or $2$.

Let $(u,u_1,u_2,u_3,u_4)$ be a $4$-geodesic of $\Ga$, and let $H=G_{u,u_1,u_2,u_3,u_4}$.
Since $\Ga$ is $(G,5)$-geodesic transitive,
we have $|H|=2^4\cdot 3$, and $H$ is transitive on $\Delta:=\Ga(u_4)\cap\Ga_5(u)$ with degree $b_4$.
For each vertex $v$ of $\Ga$, let $G_v^*$ be the kernel of $G_v$ acting on $\Ga(v)$.
Notice that $G_u^{\Ga(u)}\cong\PGL(2,13)$,  $G_u^*\cong[13^3]:\ZZ_4$ and $G_{u,u_1}\cong[13^3]:(\ZZ_{12}\times\ZZ_4)$.
Let $x\in G_u^*$ be an element of order $4$,
and let $y\in G_{u,u_1}^{\Ga(u)}$ be an element of order $12$.
Then we can write $G_{u,u_1}=[13^3]:(\la x\ra\times\la y\ra)$.
Since $|H|=2^4\cdot 3$ and $H\leq G_{u,u_1}$, we conclude that $H$ is a Hall $13'$-subgroup of $G_{u,u_1}$.
By Magma~\cite{Magma}, all Hall $13'$-subgroups of $G_{u,u_1}$ lie in the same conjugacy class.
Thus, we may choose $H=\la x\ra\times\la y\ra$.
Next, we claim that $b_4=12$.

Assume that $c_4=1$.
Then $|\Ga_4(u_0)|=14\cdot 13^3$, and so $|G_{u,u_4}|=48$ as $G_u$ is transitive on $\Ga_4(u)$.
It follows that $G_{u,u_4}=H$.
Since $\g_\Ga=8$ or $9$, we have $1\leq a_4\leq 11$.
Thus, $\Ga(u_4)\cap\Ga_4(u)$ is a union of several nontrivial orbits of $G_{u,u_4}$.
By Magma~\cite{Magma}, $G_{u,u_4}$ has one nontrivial orbit of length $1$,
nineteen orbits of length $12$, and the remaining orbits have length at least $13$.
Therefore, $a_4$ must equal to $1$, and hence $b_4=12$.

Assume that $c_4=2$.
Then $\g_\Ga=8$, $|\Ga_4(u_0)|=7\cdot 13^3$, and so $|G_{u,u_4}|=96$.
By Magma~\cite{Magma}, we have $G_{u,u_4}\cong \ZZ_4.\Sy_4$ or $\D_{24}:\ZZ_4$, and $\ZZ_4.\Sy_4$ contains no subgroups isomorphic to $H$.
Since $H\leq G_{u,u_4}$, we obtain $G_{u,u_4}\cong \D_{24}:\ZZ_4$.
Suppose that $a_4\geq1$.
Then $\Ga(u_4)\cap\Ga_4(u)$ is a union of several nontrivial orbits of $G_{u,u_4}$.
Again by Magma~\cite{Magma}, $G_{u,u_4}$ has one orbit of length $6$,
five orbits of length $12$, and the remaining orbits have length at least $13$.
If $a_4=6$, then the corresponding orbital graph has girth $3$,
which means that the subgraph $[\Ga_4(u)]$ induced by $\Ga_4(u)$ contains a cycle of length $3$.
This is impossible because $\g_\Ga=8$.
Thus, $a_4=12$, and so $b_4=0$, contradicting to $b_4\geq2$.
It follows that $a_4=0$ and so $b_4=12$.
This completes the proof of the claim.

Recall that $H=\la x\ra\times\la y \ra\cong\ZZ_4\times\ZZ_{12}$ is transitive on $\Delta$ with degree $b_4$.
By the claim, we have $b_4=12$, and so $H^\Delta=\la z\ra\cong\ZZ_{12}$.
Then $z$ must permute $\Delta$ cyclically, and so the preimage of $z$ is $(xy^3)^i$ or $y^i$ for some $i$.
Since the cycle decomposition of the element $(xy^3)^i$ is a product of two cycles of length $3$ and $4$,
we have that $\la (xy)^i\ra$ acting on $\Delta$ induces an intransitive group.
Therefore, the preimage of $z$ is $y^i$ for some $i$.
Let $K$ be the kernel of $H$ acting on $\Delta$.
Then $K=\la x\ra$.
In this case, $x$ fixes $\Delta\cup\{u_3\}$ pointwise,
and $\Ga(u_4)\setminus(\Delta\cup\{u_3\})$ contains a single vertex, which is also fixed by $x$.
Thus, $x$ fixes $\Ga(u_4)$ pointwise, which means that $x\in G_{u_4}^*$.
Recall that $x\in G_u^*$.
Thus, $x\in G_u^*\cap  G_{u_4}^*$.
By the connectivity and the $(G,5)$-geodesic transitive transitivity of $\Ga$,
$x$ fixes $V(\Ga)$ pointwise, which is clearly impossible.
This completes the proof. \qed

\begin{lemma}\label{biprimitive}
Let $\Ga$ be a connected $(G,s)$-geodesic transitive graph with girth $\g_\Ga=2s-2$ or $2s-1$,  where $G\leq\Aut(\Ga)$ and $s\geq5$.
Then $G$ cannot be primitive on $V(\Ga)$.
\end{lemma}

\proof Assume that $G$ is biprimitive on $V(\Ga)$.
Let $k$ be the valency of $\Ga$ and let $u\in V(\Ga)$.
By Lemma~\ref{sgeo-to-sarc}, $\Ga$ is $(G,s-1)$-transitive, and so it is $4$-arc transitive.
From~\cite[Theorem 1.4]{Li2001} we know that either
\begin{enumerate}[\rm (a)]
  \item $\Ga$ is the standard double cover of a vertex-primitive $s$-transitive graph $\Sigma$ (as listed in~\cite[Table 1]{Li2001}); or
  \item $\Ga\cong\Cos(G, H, HgH)$ such that $G, \soc(G)\cap H, k, s$ are as listed in~\cite[Table 2]{Li2001}.
\end{enumerate}

Assume that case (a) occurs.
By~\cite[Table 1]{Li2001}, we have $\Sig\ncong\Delta_{3,q}$, $\Delta_{4,q}$ or $\Delta_{5,q}$.
Since $\Ga$ is $(G,s)$-geodesic transitive, $\Sig$ is $s$-geodesic transitive of girth $\g_\Ga$ by Lemma~\ref{reduce}.
Moreover, $\Aut(\Sig)$ is primitive on $V(\Sig)$ in the case.
By Lemma~\ref{primitive}, such a graph $\Sig$ does exists, leading to a contradiction.

Assume that case (b) occurs. By~\cite[Table 2]{Li2001} and Proposition~\ref{geodesic-arc}, we have either
\begin{enumerate}[\rm (b.1)]
  \item $\Ga$ is geodesic transitive, and $\Ga\cong\Delta_{3,q}$, $\Delta_{4,q}$ or $\Delta_{5,q}$; or
  \item $\Ga$ is not geodesic transitive,
        and $(G,G_u,k,s-1)=(\PGL(2,p),\Sy_4,3,4)$, $(\PGammaL(2,9)$, $\Sy_4\times\Sy_2,3,5)$, or $(\M_{12}.\ZZ_2,\ZZ_3^2:\GL(2,3), 4,4)$.
\end{enumerate}
Since $\Ga$ is $(G,s)$-geodesic transitive, by Lemma~\ref{sgeo-to-sarc}, $\Ga$ is $(G,s-1)$-transitive.
If $\Ga\cong \Delta_{3,q}$, $\Delta_{4,q}$ or $\Delta_{5,q}$,
then from Table~\ref{tab:geod-sarc}, we conclude that $s=5$, $6$, $8$, respectively.
Moreover, we have $\diam(\Delta_{3,q})=3$, $\diam(\Delta_{4,q})=4$ and $\diam(\Delta_{5,q})=6$.
Therefore, $\diam(\Ga)<s$, which is impossible as $\diam(\Ga)\geq s$.
Thus, case (b.1) does not occur.

For case (b.2), we have $s=5$ or $6$.
If $s=5$, then Eq.~\eqref{eq:primitive5} also holds, and $(G,G_u,k)=(\PGL(2,p),\Sy_4,3)$ or $(\M_{12}.\ZZ_2,\ZZ_3^2:\GL(2,3),4)$.
From Eq.~\eqref{eqstab} we know that $k(k-1)^3b_4$ is a divisor of $|G_u|$, and so $b_4=1$.
If $s=6$, then Eq.~\eqref{eq:primitive6} holds and $(G,G_u,k)=(\PGammaL(2,9),\Sy_4\times\Sy_2,3)$.
By Eq.~\eqref{eqstab}, we obtain $k(k-1)^4b_5$ divides $|G_u|$, namely, $b_5=1$.
In either case, we have that $\Ga$ is geodesic transitive by Proposition~\ref{dis-tran}, leading to a contradiction.
\qed

\medskip
\noindent{\bf Proof of Theorem~\ref{Thm:quasi}.}
Combining Lemmas~\ref{reduceAS}, \ref{primitive} and~\ref{biprimitive},
we conclude that Theorem~\ref{Thm:quasi} holds. \qed

\section*{Acknowledgements}

The authors would like to thank the anonymous referee for careful reading and valuable suggestions to this paper.

\end{document}